\documentclass[12pt,a4paper]{article}
\usepackage{amsmath,amsthm,amsfonts,amssymb,bbm}

\def\RR{{\mathbb R}}
\def\CC{{\mathbb C}}
\def\NN{{\mathbb N}}
\def\ZZ{{\mathbb Z}}
\def\diff{{{\rm Diff}^+(S^1)}}
\def\mob{{\rm M\ddot{o}b}}

\def\A{{\mathcal A}}

\def\H{{\mathcal H}}

\def\M{{\mathcal M}}

\def\U{{\mathcal U}}

\newtheorem{theorem}{Theorem}[section]

\newtheorem{corollary}[theorem]{Corollary}
\newtheorem{proposition}[theorem]{Proposition}
\newtheorem{lemma}[theorem]{Lemma}

\begin{document}

\title{Restricting positive energy representations of ${\rm
Diff}^+(S^1)$ to the stabilizer of $n$ points}

\author{
{\normalsize \bf Mih\'aly Weiner\footnote
{Supported by MIUR, GNAMPA-INdAM, EU networks 
``Noncommutative Geometry'' (MRTN-CT-2006-031962)
and ``Quantum Spaces -- Noncommutative Geometry'' 
(HPRN-CT-2002-00280), and by the
``Deutsche Forschungsgemeinschaft''.}
}
\\
{\normalsize Dipartimento di Matematica,}
\\
{\normalsize Universit\`a di Roma ``Tor Vergata''}
\\
{\normalsize Via della Ricerca Scientifica, 1, I-00133, Roma, ITALY}
\\
{\small E-mail: \texttt{weiner@mat.uniroma2.it}}
}
\date{}

\maketitle

\begin{abstract}
Let $G_n \subset {\rm Diff}^+(S^1)$ be the stabilizer of $n$ given
points of $S^1$. How much information do we lose if we restrict a
positive energy representation $U^c_h$ associated to an admissible
pair $(c,h)$ of the central charge and lowest energy, to the
subgroup $G_n$? The question, and a part of the answer originate
in chiral conformal QFT.

The value of $c$ can be easily ``recovered'' from such a
restriction; the hard question concerns the value of $h$. If
$c\leq 1$, then there is no loss of information, and accordingly,
all of these restrictions are irreducible. In this work it is
shown that $U^c_{h}|_{G_n}$ is always irreducible for $n=1$ and,
if $h=0$, it is irreducible at least up to $n\leq 3$. Moreover, an
example is given for $c>2$ and certain values of $h \neq \tilde{h}$ 
such that $U^c_{h}|_{G_1}\simeq U^c_{\tilde{h}}|_{G_1}$. It is 
also concluded that for these values $U^c_{h}|_{G_n}$ cannot be 
irreducible for $n\geq 2$. For further values of $c,h$ and $n$, the 
question is left open. Nevertheless, the example already shows that in
general, local and global intertwiners in a QFT model may not be
equivalent.
\end{abstract}

\section{Introduction}

This paper concerns a purely mathematical problem regarding the
representation theory of infinite dimensional Lie groups, and it is
intended to be largely self-contained. The actual proofs, apart
from those of Prop.\! \ref{prop:c<1} and Corollary
\ref{cor:localglobal}, will not require any knowledge of chiral
conformal QFT. Nevertheless, at least in this introductory
section, we shall shortly discuss the physical motivations.

A chiral component of a conformal QFT ``lives'' on a lightline,
but it is often extended to the compactified lightline, that is,
to the circle. For several reasons, it is more convenient to study
such a theoretical model on the circle than on the lightline.
However, keeping in mind the physical motivation, one should
always clear the relation between the properties that a model has
on the lightline and the properties that it has on the circle.

For example, one may adjust the Doplicher-Haag-Roberts (DHR)
theory to describe charged sectors of models given on the circle
in the setting of Haag-Kastler nets. It is well known that a
model, when restricted to the lightline, may admit new sectors
that cannot be obtained by restrictions. These sectors are usually
called {\it solitonic}. However, so far one may have thought that
the restriction from the circle to the lightline is at least
injective: each sector restricts to a sector (i.e.\! to something
irreducible, and not to a sum of sectors), and different sectors
restrict to different sectors. In fact, under the assumption of 
{\it strong additivity}, this is indeed true. However, there are 
interesting (i.e.\! not pathological) models, in which strong additivity
fails; most prominently, the Virasoro net with central charge
$c>1$. (The Virasoro nets are fundamental, because each chiral
conformal model contains a Virasoro net as a subnet in an
irreducible way.)

As it is known, see e.g.\! the book \cite{kac}, for certain values
of the central charge $c$ and lowest energy $h$, there exists a
unitary lowest energy representation of the Virasoro algebra. By
\cite{GoWa}, each of these representations gives rise to a
projective unitary representation $U^c_h$ of the group $\diff$,
that is, of the group of orientation-preserving smooth
diffeomorphisms of the circle. These representations are all
irreducible, and every positive energy irreducible representation
of $\diff$ is equivalent with $U^c_h$ for a certain admissible
pair $(c,h)$. Moreover, two of these representations are
equivalent if and only if both their central charges and their
lowest energies coincide.

A representation $U^c_h$ with lowest energy $h=0$ gives rise to a
conformal net (in its vacuum representation) on the circle. This
is the so-called Virasoro net at central charge $c$, and it is
denoted by $\A_{{\rm Vir}_c}$. Every charged sector of $\A_{{\rm
Vir}_c}$ arises from a positive energy representation of $\diff$
with (the same) central charge $c$. Two charged sectors are
equivalent if and only if they arise from equivalent positive
energy representations of $\diff$.

Viewing the circle as the compactified lightline, i.e.\! the
lightline together with the ``infinite'' point, one has that a
diffeomorphism of the circle restricts to a diffeomorphism of the
lightline if and only if it stabilizes the chosen infinite point.
So by what was roughly explained, we are motivated to ask the
following questions. Let $G_n \subset \diff$ be the stabilizer
subgroup of $n$ given points of $S^1$. Then

- is the restriction of $U^c_h$ to $G_1$ irreducible?

- for what values of $(c,h)$ and $(\tilde{c},\tilde{h})$ we have
$U^c_h|_{G_1} \simeq U^{\tilde{c}}_{\tilde{h}}|_{G_1}$?

\noindent Actually, with $G_1$ replaced with $G_n$, there are
reasons to consider these questions not only for $n=1$, but in
general. (Note that though the actual elements of $G_n$ depend on
the choice of the $n$ points, different choices result in
conjugate subgroups: thus all of these questions are well-posed.)
In fact, the (possible) irreducibility of $U^c_h|_{G_n}$ for
$h=0$, is directly related to the (possible) $n$-regularity of
$\A_{{\rm Vir}_c}$. (See \cite{GLW} for more on the notion of
$n$-regularity.) The other reason is the relation between the
answers regarding different values of $n$. Of course, we have some
trivial relations, since $G_n$ may be considered to be a subgroup
of $G_m$ whenever $n\geq m$. However, as it will be proved at
Prop.\! \ref{prop:geometrical}, we have the further relation:
\begin{equation}\nonumber
U^c_h|_{G_{n+1}} \;\;{\rm is\;\;irred.} \;\;\;\;\Rightarrow \;\;\;\;
U^c_h|_{G_n}\simeq U^{\tilde{c}}_{\tilde{h}}|_{G_n}
\;\;{\rm if  \;\; and\;\;only \;\;if}\;\;(c,h)=
(\tilde{c},\tilde{h}).
\end{equation}

As it was mentioned, both the questions, and a part of their answers
originate in chiral conformal QFT. For example, {\it Haag-duality} is
known \cite{BGL,FrG} to hold in the vacuum sector of any chiral conformal
net on the circle. This could be used to conclude that $U^c_0|_{G_2}$ is
irreducible for all values of the central charge $c$. However, we shall
not enter into details of this argument, because in any case we shall
prove some stronger statements regarding irreducibility. In particular, by
considering the problem at the Lie algebra level, it will be shown that
the the representation $U^c_h|_{G_n}$ for $n=1$, is always irreducible
(Corollary \ref{cor:irred}), and for $h=0$, we have irreducibility at
least up to $n\leq 3$ (Prop.\! \ref{prop:h=0}).

Apart from general statements regarding conformal nets, by now we
have a detailed knowledge, in particular, of Virasoro nets. For
example, it is known, that for $c\leq 1$ they are strongly
additive, \cite{kl,xu03}. This permits us to conclude (Prop.\!
\ref{prop:c<1}), that whenever $c\leq 1$, the representation
$U^c_h|_{G_{n}}$ is irreducible for any positive integer $n$, and
accordingly, $U^c_h|_{G_n}\simeq U^{\tilde{c}}_{\tilde{h}}|_{G_n}$
if and only if $(c,h) = (\tilde{c},\tilde{h})$.

Thus for $c\leq 1$, our questions are answered. Let us discuss now the
region $c>1$. It is easy to show, that --- in general ---
$U^c_h|_{G_n}\simeq U^{\tilde{c}}_{\tilde{h}}|_{G_n}$ implies
$c=\tilde{c}$ (Corollary \ref{cor:c=c}). Hence the real problem is to
``recover'' the value of the lowest energy.

The main result of this paper is an example, showing that already
for $n=1$, the value of the lowest energy cannot be always
determined by the restriction, since in particular for
$h=\frac{c}{32}, {\tilde{h}}=\frac{c}{32} + \frac{1}{2}$ and $c>2$
we have $U^c_h|_{G_1}\simeq U^{\tilde{c}}_{\tilde{h}}|_{G_1}$. It
follows (Corollary \ref{cor:localglobal}), that the Virasoro net
with $c>2$ is an example for a net in which local and global
intertwiners are not equivalent.

The values here exhibited may have more to do with the actual
construction, than with the problem itself. (For example, it could
turn out that when $c>1$, all of the representation $U^c_h|_{G_1}$
with $h$ varying over the positive numbers, are equivalent.) The
method of showing equivalence is obtained by a combined use of two
known tricks:

- the realization (appearing e.g.\! in \cite{BS-M}) of the
lightline-restriction of the Virasoro net at $c>1$ as a subnet of the so
called $U(1)$-current,

- the observation (appearing e.g.\! in \cite{LX}) that for any
$k=1,2,...$, the map $l_n\mapsto \frac{1}{k}l_{kn} +
\frac{C}{24}(k-\frac{1}{k})\delta_{n,0}$ gives an endomorphism of
the Virasoro algebra.

\section{Preliminaries}

The Virasoro algebra (${\rm Vir}$) is spanned by the elements $\{l_n : n
\in \ZZ \}$ together with the central element $C$ obeying the commutation
relations
\begin{eqnarray}
\nonumber
[l_n,l_m] &=& (n-m)l_{n+m} + \frac{C}{12} (n^3-n) \delta\!_{-n,m}
\\ {[l_n,C]} & = & 0.
\end{eqnarray}
For a representation $\pi$ of ${\rm Vir}$ on a complex vector space
$V$, set $L_n\equiv\pi(l_n)$. An eigenvalue of $L_0$ is usually
referred as a value of the {\bf energy}, and the corresponding eigenspace
as the {\bf energy level} associated to that value. If $L_0\Phi=\lambda
\Phi$, then by a use of the commutation relations
$L_0(L_n\Phi)=(\lambda-n)\Phi$, i.e.\! the operator $L_n$ decreases the
value of the energy by $n$.
We say that $\pi$, with representation space $V\neq 0$, is a {\bf lowest 
energy representation} with {\bf central charge} $c\in \CC$ and {\bf 
lowest energy} $h\in \CC$, iff

(1) $h$ is an eigenvalue of $L_0$, and if ${\rm Re}(s)<{\rm Re}(h)$ then
$s$ is not an eigenvalue of $L_0$ (i.e.\! $h$ is the ``lowest energy''),

(2) $\pi(C)=c\mathbbm 1$,

(3) $V$ is spanned by the orbit (under $\pi$) of a 
vector of ${\rm Ker}(L_0-h\mathbbm 1)$.

\noindent In this case ${\rm Ker}(L_0-h\mathbbm 1)$ is
one-dimensional, so up to a multiplicative constant there exists a
unique {\bf lowest energy vector} $\Psi^c_h$ and $L_n\Psi^c_h=0$ for all
$n>0$. Moreover, $V=\oplus_{n=0}^\infty V_{(h+n)}$ where
$V_{(h+n)}={\rm Ker}(L_0-(h+n)\mathbbm 1)$ for $n=0,1,...$ and
actually the dimension of $V_{(h+n)}$ is smaller than or equal to
the number of partitions of $n$, as in fact
\begin{equation}
V_{(h+n)}={\rm Span}\{L_{-n_1}...L_{-n_j}\Psi^c_h|\, j\in\NN,\,
n_1\geq...\geq n_j >0, \sum_{l=0}^jn_l=n\}
\end{equation}
where $j=0$ means that no operator is applied to $\Psi^c_h$.

A {\bf unitary representation} of the Virasoro algebra
is a representation $\pi$ of ${\rm Vir}$ on
a complex vector space $V$ endowed with a (skew symmetric,
positive definite) scalar product $\langle\cdot,\!\cdot\rangle$
satisfying the condition
\begin{equation}
\langle\pi(l_n)\Phi_1,\Phi_2\rangle =
\langle\Phi_1,\pi(l_{-n})\Phi_2\rangle \;\;\;\;(\Phi_1,\Phi_2\in
V,\;n\in\ZZ),
\end{equation}
or in short, that $\pi(l_n)^+\equiv \pi(l_n)^*|_{V}=\pi(l_{-n})$. 
(We use the symbol ``$^+$'', keeping ``$^*$'' exclusively for the adjoint
defined in the von Neumann sense on a Hilbert space.) Note that
the formula $\theta(l_n)=l_{-n}$ defines a unique antilinear
involution with the property that $[\theta(x),\theta(y)]=\theta([y,x])$,
and that unitaritity means that $\pi(x)^+=\pi(\theta(x))$ for every
$x\in {\rm Vir}$.

A pair $(c,h)$ is called {\bf admissible}, if there exists a unitary
lowest energy representation with central charge $c$ and lowest
energy $h$. If $(c,h)$ is admissible, then up to equivalence there
exists a unique unitary lowest energy representation with central
charge $c$ and lowest energy $h$. In this paper this unique
representation will be denoted by $\pi^c_h$, the corresponding
representation space by $V^c_h$, and the (up-to-phase unique)
normalized lowest energy vector by $\Psi^c_h$. As is known,
this representation is irreducible (in the algebraic sense) and
two such representations are equivalent (in the algebraic sense)
if and only if their central charges, as well as their lowest
energies, coincide.

Of course $(c,h)=(0,0)$ is an admissible pair and the
corresponding representation is trivial, but a pair $(c,h)\neq
(0,0)$, as is known (see e.g.\! the book \cite{kac} for further
explanations), is admissible if and only if it belongs to either
to the {\it continuous part} $[1,\infty)\times [0,\infty)$ or to
the {\it discrete part} $\{(c(m),h_{p,q}(m))|m\in\NN,\,
p=1,...,m+1;\, q=1,...,p\}$ where
\begin{equation}
c(m)=1-\frac{6}{(m+2)(m+3)},\;\;\;
h_{p,q}(m)=\frac{((m+3)p-(m+2)q)^2-1}{4(m+2)(m+3)}.
\end{equation}

Let us see now what all this has to do with the so-called positive
energy representations of $\diff$, where by the symbol ``$\diff$''
we mean the group of orientation preserving (smooth)
diffeomorphisms of the unit circle $S^1 \equiv \{z\in \CC|\,
\|z\|=1\}$. We shall always consider $\diff$ as a continuous group
with the usual $C^\infty$ topology.

We shall often think of a smooth function $f\in C^\infty(S^1,\RR)$ as the
vector field symbolically written as $z=e^{i\theta}\mapsto
f(e^{i\vartheta})\frac{d}{d\vartheta}$. The corresponding one-parameter
group of diffeomorphisms will be denoted by $t\mapsto {\rm Exp}(tf)$.

We shall denote by $\U(\H)$ the group of unitary operators of a
Hilbert space $\H$. A {\bf projective unitary operator} on $\H$ is
an element of the quotient group $\U(\H)/\{z{\mathbbm 1}|z\in
S^1\}$. A (strongly continuous) projective representation of a
(continuous) group $G$ is a (strongly continuous) homomorphism
from $G$ to $\U(\H)/\{z{\mathbbm 1}|z\in S^1\}$.

We shall often think of a projective unitary operator $Z$ as a
unitary operator. Although there are more than one way of fixing
phases, note that expressions like ${\rm Ad}(Z)$ or $Z \in \M$ for
a von Neumann algebra $\M \subset {\rm B}(\H)$ are unambiguous.
Note also that the self-adjoint generator of a one-parameter group
of strongly continuous {\it projective} unitaries $t \mapsto Z(t)$
is well defined up to a real additive constant: there exists a
self-adjoint operator $A$ such that ${\rm Ad}(Z(t))= {\rm
Ad}(e^{iAt})$ for all $t \in \RR$, and if $A'$ is another
self-adjoint with the same property then $A'=A+r\mathbbm 1$ for
some $r \in \RR$.

Let now $(c,h)$ be an admissible pair for ${\rm Vir}$, and denote
by $\H^c_h$ the Hilbert space obtained by the completion of the
representation space $V^c_h$ of $\pi^c_h$. The operator
$L_n=\pi^c_h(l_n)$ may be viewed as a densely defined operator on
this space. By the unitarity of $\pi^c_h$, we have that
$L_n^*\supset L_{-n}$ (i.e.\! $L_n^*$ is an extension of $L_{-n}$)
and hence $L_n$ is closable.

If $f:S^1 \rightarrow \CC$ is a smooth function with Fourier
coefficients
\begin{equation}
\hat{f}_n\equiv \frac{1}{2\pi}\int_0^{2\pi}
e^{-in\theta}f(e^{i\theta})d\theta \;\;\;\;(n\in\ZZ),
\end{equation}
then the sum $\sum_{n\in\ZZ} \hat{f}_n L_n$ is strongly convergent
on $V^c_h$, and the operator given by the sum is closable.
Denoting by $T^c_h(f)$ the corresponding closed operator, by a use
of Nelson's commutator theorem \cite[Prop.\! 2]{Ne}, one has that 
$T^c_h(f)^* =
T^c_h(\overline{f})$ and so in particular that $T^c_h(f)$ is
self-adjoint whenever $f$ is a real function. By the main result
of \cite{GoWa}, there exists a unique projective unitary
representation $U^c_h$ of $\diff$ on $\H$ such that
\begin{equation}
U^c_h({\rm Exp}(f)) = e^{iT^c_h(f)}
\end{equation}
for every $f\in C^\infty(S^1,\RR)$. This representation is
strongly continuous, and moreover, it is irreducible. Note that
this latter property does not follow immediately from the fact
that $\pi^c_h$ is irreducible (in the algebraic sense). For
example, $U^c_h$ could have a nontrivial invariant closed subspace
which has a trivial intersection with the dense subspace $V^c_h$;
see also the related remark after Prop.\! \ref{prop:equivalence}.
However, this is not so. Indeed, if $W$ is bounded operator in the
commutant of $U^c_h$, then, in particular $W\overline{L}_0\subset
\overline{L}_0W$ and hence $W$ preserves each eigenspace of
$\overline{L}_0$. But by assumption one can forme a complete
orthonormed system consisting of eigenvectors of $L_0$, and so the
eigenspaces of $\overline{L}_0$ are exactly the eigenspaces of
$L_0$. Thus $W$ preserves each energy space and so also the dense
subspace $V^c_h$. Similar arguments show that $U^c_h\equiv
U^{\tilde{c}}_{\tilde{h}}$ if and only if $c=\tilde{c}$ and
$h=\tilde{h}$.

A {\bf positive energy representation} $U$ of $\diff$ on $\H$ is a
strongly continuous homomorphism from $\diff$ to
$\U(\H)/\{z{\mathbbm 1}|z\in S^1\}$ such that the self-adjoint
generator of the anticlockwise rotations is bounded from below.
(Note that although the generator is defined only up to a real
additive constant, the fact whether it is bounded from below is
unambiguous.)

Diffeomorphisms of $S^{1}$ of the form $z \mapsto
\frac{az+b}{\overline{b}z+\overline{a}}$ with $a,b\in \CC$,
$|a|^2-|b|^2=1$ are called {\bf M\"obius-transformations}. The
subgroup $\mob\subset \diff$ formed by these transformations is
isomorphic to ${\rm PSL}(2,\RR)$, and it is generated by the (real
combinations of the complex) vector fields $z\mapsto z^{\pm 1}$
and $z\mapsto 1$. Note that in the representation $U^c_h$, the
three listed complex vector fields correspond to the three
operators $L_{\pm 1}$ and $L_0$.

A strongly continuous projective representation of $\mob$ always
lifts to a unique strongly continuous unitary representation of
the universal covering group $\widetilde{\mob}\equiv
\widetilde{{\rm PSL}(2,\RR)}$ of $\mob$. Through restriction and
this lifting, one may fix the additive constant in the definition
of the self-adjoint generator of anticlockwise rotations and
define the {\bf conformal Hamiltonian} $L_0$ of a strongly
continuous projective representation of $\diff$. As is well
known, $L_0$ is bounded from below (i.e.\! the representation is
of positive energy type) if and only if $L_0$ is actually bounded
by $0$. Moreover, each irreducible positive energy representation
of $\diff$ is equivalent to $U^c_h$ for a certain admissible pair
(c,h); see \cite[Theorem A.2]{Car}.

\section{The passage to the Lie algebra level}

One could view $G_n$ as a Lie subgroup of $\diff$, with the
corresponding Lie algebra consisting of those vector fields that
vanish at the given $n$ points. Without any loss of generality,
let us assume that the given points of the unit circle are
$e^{i\frac{2\pi}{n}k}$ for $k=1,...,n.$ Then the mentioned
(complexified) Lie subalgebra can be identified with the set of
functions ${\mathfrak G}_n \equiv \{f\in C^\infty(S^1,\CC):
f(e^{\frac{2\pi}{n}k})=0\;{\rm for}\; k=1,...,n\}.$ To find a
suitable base, consider the function defined by the formula
\begin{equation}
e_{j,r}(z)\equiv z^j-z^{r+j}= z^j(1-z^{r})
\end{equation}
where $r\in \ZZ$ and $j\in \NN$. Then $\{e_{j,rn}: r\in\ZZ,
j=0,...,n-1\}$ is a set of linearly independent elements of
${\mathfrak G}_n$ whose span is dense in ${\mathfrak G}_n$, where
the latter is considered with the usual $C^\infty$ topology.

Omitting the indices of central charge and lowest energy, we have
that $T(e_{j,r})=\overline{\pi(l_j-l_{r+j})}$. So let us set
\begin{equation}
k_{j,r} \equiv l_j - l_{r+j}.
\end{equation}
We shall often use $k_{0,r}$. To shorten formulae, we shall set
$k_r\equiv k_{0,r}$. By direct calculation we find that
\begin{equation}\label{[k,k]}
[k_r,k_m] = r k_r - m k_m - (r-m) k_{r+m} + \frac{C}{12}
(r^3-r) \delta\!_{-r,m},
\end{equation}
implying that the elements $\{k_r : r \in\ZZ\}$ together with the
central element $C$ span a Lie subalgebra of the Virasoro algebra,
which we shall denote by $\mathfrak K$. In fact, by a similar
straightforward calculation one has that
\begin{equation}
\mathfrak K_n \equiv {\rm Span}(\{k_{j,rn}: r\in \ZZ,
j=0,...,n-1\}\cup \{C\})
\end{equation}
is a Lie subalgebra of ${\rm Vir}$ for any positive integer $n$. (Note
that for $n=1$ we get back $\mathfrak K$, i.e.\! ${\mathfrak K}=\mathfrak
K_1$.) Intuitively, viewing the Virasoro algebra from the point of view of
vector fields on the circle, ${\mathfrak K}_n$ corresponds to the algebra
of Laurent-polynomial (polynomial in $z$ and $z^{-1}$)
vector fields, that are zero at the chosen $n$ points of $S^1$.
In what follows, and throughout the rest of this paper, for a densily 
defined operator $A$ we shall denote its closure by $\overline{A}$.

\begin{proposition}
\label{prop:equivalence}
Let $(c,h)$ and $(\tilde{c},\tilde{h})$ be two admissible pairs for the
Virasoro algebra, and assume that $U^c_h|_{G_{n}}\simeq U^{\tilde
c}_{\tilde h}|_{G_{n}}$. Then there exists a unitary operator $V: \H^c_h
\to \H^{\tilde c}_{\tilde h}$ and a linear functional
$\phi:{\mathfrak K}_n \to \CC$ with ${\rm Ker}(\phi)\supset
[{\mathfrak K}_n, {\mathfrak K}_n]$
such that for all $x\in {\mathfrak K}_n$ we have
$$
V \overline{\pi^c_h(x)}V^* =
\overline{\pi^{\tilde{c}}_{\tilde{h}} (x)} +\phi(x)\mathbbm 1.
$$
\end{proposition}
\begin{proof}
If the real vector field $f$ belongs to ${\mathfrak G}_n$, then
${\rm Exp}(tf)\in G_n$ for every $t\in \RR$. Using that both the
real part and the imaginary part of $x$ is in ${\mathfrak G}_n$,
the fact that the finite energy vectors form a core for all
operators of the form $T(f)$, and some standard arguments, one can
easily show that if the two representations are equivalent, then
there exists a unitary $V$ such that $V \overline{\pi^c_h(x)} V^*
= \overline{\pi^{\tilde{c}}_{\tilde{h}} (x)}\, +$ an additive
constant, that may depend (linearly) on $x$; say $\phi(x)\mathbbm
1$. (Recall that we are dealing with {\it projective}
representations.) As $\pi$ is a Lie algebra representation, we
have that $[\pi(x),\pi(y)]=\pi([x,y])$. Actually, as
$\overline{\pi(x)}=\pi(\theta(x))^*$, we have that
\begin{eqnarray}\nonumber
\,[\overline{\pi(x)},\overline{\pi(y)}] &=&
[\pi(\theta(x))^*,\pi(\theta(y))^*] \\
&\subset&
[\pi(\theta(y)),\pi(\theta(x))]^* =
\pi(\theta([x,y]))^* =\overline{\pi([x,y])}.
\end{eqnarray}
Thus it follows that
\begin{eqnarray}\nonumber
V\pi^c_h([x,y])V^* &=&
V[\pi^c_h(x),\pi^c_h(y)]V^*=
[V\pi^c_h(x)V^*,V\pi^c_h(y)V^*] \\
\nonumber
&\subset&
[(\overline{\pi^{\tilde{c}}_{\tilde{h}}(x)}+\phi(x)\mathbbm 1),
(\overline{\pi^{\tilde{c}}_{\tilde{h}}(y)}+\phi(y)\mathbbm 1)]
=[\overline{\pi^{\tilde{c}}_{\tilde{h}}(x)},
\overline{\pi^{\tilde{c}}_{\tilde{h}}(y)}]
\\
&\subset& \overline{\pi^{\tilde{c}}_{\tilde{h}}([x,y])}
\end{eqnarray}
and hence $\phi([x,y])=0$, which concludes our proof.
\end{proof}
{\it Remark}. Note that even if $\phi=0$, the unitary operator $V$
appearing in the above proposition does {\it not}\, necessarily
make an equivalence between $\pi^c_h|_{{\mathfrak K}_n}$ and
$\pi^{\tilde{c}}_{\tilde{h}}|_{{\mathfrak K}_n}$, since it may not
take the dense subspace $V^c_h$ into $V^{\tilde{c}}_{\tilde{h}}$.

This possibility is not something which is specific to infinite
dimensional Lie groups. In fact, consider two of the unitary
lowest energy irreducible representations (with lowest energy
different from zero), say $\eta_1$ and $\eta_2$ of the Lie algebra
${\mathfrak{sl}}(2,\RR)$. Moreover, consider the base $e_+,e_-$
and $h$ (in the complexified) Lie algebra satisfying the usual
commutation relations $[h,e_\pm]= \mp e_\pm$ and $[e_-,e_+]=2h$.
The two elements $t=2h-(e_-+e_+)$ and $s=i(e_--e_+)$ span a
two-dimensional Lie subalgebra of ${\mathfrak{sl}}(2,\RR)$, and it
is easy to prove, that if $\eta_1$ and $\eta_2$ are inequivalent,
then also their restrictions to this subalgebra are inequivalent
(in the algebraic sense). However, the corresponding
representations of the corresponding Lie sub{\it{groups}} are in
fact equivalent; see for example the remarks in the proof of
\cite[Theorem 2.1]{GLW}.
\bigskip

\begin{proposition}\label{prop:[K,K]}
$C,k_{rn}\in[{\mathfrak K}_n,{\mathfrak K}_n]$ for every $r\in \ZZ$ and
positive integer $n$. In particular, $[{\mathfrak K},{\mathfrak K}]
={\mathfrak K}$.
\end{proposition}
\begin{proof}
Let $\phi:{\mathfrak K}_n\to \CC$ be a linear functional such that
${\rm Ker}(\phi)=[{\mathfrak K}_n,{\mathfrak K}_n]$.
Our aim is to show that $\phi(C)=\phi(k_{rn})=0$.
To shorten notations, we shall set $\phi_r\equiv \phi(k_{rn})$. Then
by equation (\ref{[k,k]}) one finds that
\begin{equation}\label{s_r}
r \phi_r - m \phi_m - (r-m)\phi_{r+m} +
\frac{\phi(C)}{12}(n^2r^3-r)\delta_{-r,m}
=0
\end{equation}
for all $r,m\in \ZZ$. Let us now analyze the above relation (together with
the fact that $\phi_0=\phi(k_0)=\phi(0)=0)$. If $r>1$ and $m=1$, then by
substituting into (\ref{s_r}) we obtain the recursive relation $r\phi_r
+\phi_1 - (r-1)\phi_{r-1}=0$.
Resolving the recursive relation we get that for $r>1$ we have
\begin{equation}\label{sol:rec}
\phi_r=(r-1)(\phi_2-\phi_1)+\phi_1.
\end{equation}
Similarly, letting $r<-1$ and $m=-1$ and resolving the resulting
recursive relation we get that $\phi_r$ is a (possibly different)
first order polynomial of $r$ for the region $r<-1$, too. Then
letting $m=-r$ and using that $\phi_0=0$, we find by substitution
that
\begin{equation}
r\phi_r + r \phi_{-r} + \frac{\phi(C)}{12}(n^2r^3-r)=0.
\end{equation}
The expression on the left-hand side --- by what was just
explained --- for the region $r>1$, is a polynomial of $r$. Thus
each coefficient of this polynomial must be zero, and hence, by
what was just obtained about degrees, we find that $\phi(C)=0$
which then by the above equation further implies that
$\phi_r=-\phi_{-r}$ for every $r\in \ZZ$. Moreover, returning to
(\ref{s_r}), we have that
\begin{equation}
r\phi_r-m\phi_m-(r-m)\phi_{r+m}=0
\end{equation}
and also, by exchanging $m$ with $-m$ and using
that $\phi_{-m}=-\phi_m$, we have that
\begin{equation}
r\phi_r-m\phi_m-(r+m)\phi_{r-m}=0.
\end{equation}
Taking the difference of these two equations and setting $m=r-1$,
we find that $\phi_{2r-1}=(2r-1)\phi_1$. Restricting our attention
to the region $r>1$, and confronting what we have just obtained
with (\ref{sol:rec}), we get that $\phi_1=\phi_2=0$ and hence
again by (\ref{sol:rec}) that $\phi_r=0$ for all $r\geq 0$ and so
actually for all $r\in \ZZ$, which concludes our proof.
\end{proof}

\begin{corollary}\label{cor:c=c}
Let $(c,h)$ and $(\tilde{c},\tilde{h})$ be two admissible pairs for the
Virasoro algebra, and assume that $U^c_h|_{G_{n}}\simeq U^{\tilde
c}_{\tilde h}|_{G_{n}}$. Then $c=\tilde{c}$.
\end{corollary}
\begin{proof}It follows trivially from Prop.\! \ref{prop:equivalence} 
and \ref{prop:[K,K]}.
\end{proof}
We shall now formulate a useful condition of irreducibility. Fix
an admissible pair $(c,h)$ of the Virasoro algebra, and a positive
integer $n$. Recall that we have denoted by $\Psi^c_h$ the (up to
phase) unique normalized lowest energy vector of the
representation $\pi^c_h$. It is clear that the subset of
${\mathfrak K}_n$
\begin{equation}
{\mathfrak O}^c_{h,n}\equiv \{x\in {\mathfrak K}_n:\,
\pi(x)\Psi^c_h = \lambda_x\Psi^c_h \text{ for some
}\lambda_x\in\CC\}
\end{equation}
is in fact a Lie subalgebra. Note that $\theta({\mathfrak
K}_n)={\mathfrak K}_n$, but $\theta({\mathfrak O}^c_{h,n})\neq
{\mathfrak O}^c_{h,n}$.
\begin{proposition}
\label{prop:irred}
Suppose that $V^c_h$ is spanned by vectors of the form
$$\pi(\theta(x_1))...\pi(\theta(x_j))\Psi^c_h,$$
where $x_1,...,x_j\in {\mathfrak O}^c_{h,n}$ and $j\in \NN$ (with
$j=0$ meaning the vector $\Psi^c_h$ itself). Then $U^c_h|_{G_n}$
is irreducible.
\end{proposition}
\begin{proof}
Simple arguments (similar to those appearing in the proof of
Prop.\! \ref{prop:equivalence}) show, that if $V$ is a unitary
operator commuting with $U^c_h(G_n)$, then for all $x\in
{\mathfrak K}_n$ we have
$V\overline{\pi^c_h(x)}=\overline{\pi^c_h(x)}V$. Let $B\equiv V
-\langle \Psi^c_h,V\Psi^c_h\rangle \mathbbm 1$; then for every
$j\in \NN$ and $x_1,...,x_j\in {\mathfrak O}^c_{h,n}$ we have
\begin{eqnarray}\nonumber
\langle \pi^c_h(\theta(x_1))...\pi^c_h(\theta(x_j))\Psi^c_h,
B\Psi^c_h \rangle &=& \langle \pi(x_1)^*...\pi(x_j)^*\Psi^c_h,
B\Psi^c_h \rangle \\
\nonumber &=& \langle \Psi^c_h,
B\pi^c_h(x_j)...\pi^c_h(x_1)\Psi^c_h \rangle
\\
&=&\text{ multiple of }\langle B\Psi^c_h,\Psi^c_h\rangle =0
\end{eqnarray}
and hence by the condition of the proposition $B\Psi^c_h=0$. In
turn this implies that
$B\pi^c_h(\theta(x_1))...\pi^c_h(\theta(x_j))\Psi^c_h=
\pi^c_h(\theta(x_1))...\pi^c_h(\theta(x_j))B\Psi^c_h=0$ and hence
that $B=0$; i.e.\! that
$V=\langle\Psi^c_h,V\Psi^c_h\rangle\mathbbm 1$.
\end{proof}
In order to use the above proposition, let us fix a certain admissible
value of $c$ and $h$.
To simplify notations, we shall set $L_n \equiv \pi^c_h(l_n)$
and $K_n\equiv \pi^c_h(k_n)=L_0-L_n$. Note that $K_0 = 0$ and that
$\pi^c_h(C) = c \mathbbm 1$.
\begin{lemma}
\label{vermalike-module}
The vectors of the form
$$
K_{-n_1} K_{-n_2} \ldots K_{-n_k} \Psi^c_h
$$
where $k\in \NN, n_j \in \NN^+ (j=1\ldots k)$ and $n_1 \geq n_2 \geq
\ldots \geq n_k$ (and where $k=0$ means the vector $\Psi^c_h$ in itself,
without any operator acting on it) span the representation space $V^c_h$.
\end{lemma}
\begin{proof}
The statement with ``$K$'' everywhere replaced by ``$L$''
is true by definition. On the other hand,
\begin{equation}
L_{-n_1} \Psi^c_h = (h\mathbbm 1 - L_0 + L_{-n_1})\Psi^c_h =
-K_{-n_1}\Psi^c_h + h\Psi^c_h
\end{equation}
Similarly,
\begin{eqnarray}
\nonumber
L_{-n_1} L_{-n_2}\Psi^c_h &=& ((h+n_2)\mathbbm 1 - L_0 +
L_{-n_1})(h\mathbbm 1
-L_0+L_{-n_2})\Psi^c_h \\ \nonumber
&=&((h+n_2)\mathbbm 1 - K_{-n_1})(h\mathbbm 1
-K_{-n_2})\Psi \\
&=& K_{-n_1}K_{-n_2}\Phi_h - (h+n_2)K_{-n_2}\Psi^c_h + h\Psi^c_h,
\end{eqnarray}
and it is not too difficult to generalize the above argument, by
induction, to show that $L_{-n_1} L_{-n_2}\ldots L_{-n_k}\Psi^c_h$
is a linear combination of vectors of the discussed form.
\end{proof}
The lowest energy vector $\Psi^c_h$,
though (in general) it is not annihilated by the operators $K_n$ 
$(n>0)$, is still a common eigenvector for them: 
\begin{equation}
\label{lowestenergy-vector}
\forall n>0: \, K_n \Psi^c_h = h \Psi^c_h.
\end{equation}
Hence $\theta(k_{-n})=k_n\in {\mathfrak O}^c_{h,1}$ and thus by
Lemma \ref{vermalike-module} and Prop.\! \ref{prop:irred} we can draw
the following conclusion.
\begin{corollary}\label{cor:irred}
Let $(c,h)$ be any admissible pair.
If $n=1$ then $U^c_h|_{G_n}$
is irreducible.
\end{corollary}
By Prop.\! \ref{prop:equivalence} and \ref{prop:[K,K]}, if $V$
is a unitary operator making an equivalence between $U^c_h|_{G_1}$
and $U^{\tilde{c}}_{\tilde{h}}|_{G_1}$, then it also makes an
equivalence between $\overline{\pi^c_h}|_{\mathfrak K}$ and
$\overline{\pi^{\tilde{c}}_{\tilde{h}}}|_{\mathfrak K}$. To show the
converse, one needs to overcome the following difficulty: we do not know, 
whether the subgroup of $G_n$ generated by the exponentials is dense in 
$G_n$. In what follows we shall denote this subgroup by $\tilde{G}_n$.

\begin{proposition}\label{intertwiners}
Let $B$ be a bounded operator from $\H^c_h$ to
$\H^{\tilde{c}}_{\tilde{h}}$. Then $B$ intertwines $U^c_h|_{G_1}$
with $U^{\tilde{c}}_{\tilde{h}}|_{G_1}$ if and only if it
intertwines $U^c_h|_{\tilde{G}_1}$ with
$U^{\tilde{c}}_{\tilde{h}}|_{\tilde{G}_1}$.
\end{proposition}
\begin{proof}
Clearly, we have never used in our proof of irreducibility the whole group
$G_1$, but only the subgroup $\tilde{G}_1$ generated by the exponentials.
Hence we have that also $U^{\tilde{c}}_{\tilde{h}}|_{\tilde{G}_1}$ and
$U^{\tilde{c}}_{\tilde{h}}|_{\tilde{G}_1}$ are irreducible
representations.

If $B$ intertwines $U^c_h|_{G_1}$ with
$U^{\tilde{c}}_{\tilde{h}}|_{G_1}$ then of course it also
intertwines $U^c_h|_{\tilde{G}_1}$ with
$U^{\tilde{c}}_{\tilde{h}}|_{\tilde{G}_1}$. So let us assume that
$B$ is an intertwiner of the latter two. Then, since these
representations are unitary and irreducible, it follows that $B$
is a multiple of unitary operator. 

So we may assume that $B$ is unitary. Then, using the intertwining
property and the fact that the conjugate of an exponential in
$G_1$ is still an exponential, it is easy to show that ${\rm
Ad}(BU^c_h(g)B^*)(U^{\tilde{c}}_{\tilde{h}}(\tilde{g}))=
U^{\tilde{c}}_{\tilde{h}}(g\tilde{g}g^{-1})$ for all $g\in G_1$
and $\tilde{g}\in\tilde{G}_1$. Thus by the irreducibility of
$U^{\tilde{c}}_{\tilde{h}}|_{\tilde{G}_1}$ it follows that ${\rm
Ad}(BU^c_h(g)B^*)={\rm Ad}(U^{\tilde{c}}_{\tilde{h}}(g))$ and so
that in the projective sense
$BU^c_h(g)B^*=U^{\tilde{c}}_{\tilde{h}}(g)$, which finishes our
proof.
\end{proof}

\begin{corollary}
\label{algebra->group}
Let $(c,h)$ and $(\tilde{c},\tilde{h})$ be two admissible pairs for the
Virasoro algebra. Then $U^c_h|_{G_{1}}\simeq U^{\tilde
c}_{\tilde h}|_{G_{1}}$ if and only if there exists a unitary operator
$V$ such that $V\overline{\pi^c_h(x)}V^* =
\overline{\pi^{\tilde{c}}_{\tilde{h}}(x)}$ for all $x\in {\mathfrak K}$.
\end{corollary}
\begin{proof}
The ``only if'' part follows from  Prop.\! \ref{prop:equivalence} and 
\ref{prop:[K,K]}. As for the ``if'' part: it is clear, that if the two 
representations are equivalent at the Lie algebra level, then they are 
also equivalent on the subgroup generated by the exponentials. Hence the 
``if'' part follows directly from the previous proposition.
\end{proof}

\section{Further observations}

\begin{proposition}
\label{prop:geometrical} Let $(c,h)$ and $(\tilde{c},\tilde{h})$
be two admissible pairs, $n$ a positive integer, and suppose that
$U^c_h|_{G_{n+1}}$ is irreducible. Then $U^c_h|_{G_n}\simeq
U^{\tilde{c}}_{\tilde{h}}|_{G_n}$ if and only if
$(c,h)=(\tilde{c},\tilde{h})$.
\end{proposition}
\begin{proof}
The ``if'' part is trivial; we only need to prove the ``only if''
part. So suppose the two representations of $G_n$ in question are
equivalent. In fact, assume that they actually coincide. (Clearly,
we can safely do so.) So we shall fix $n$ (different) points
$p_1,...,p_n$ on the circle, we shall think of $G_n$ as the their
stabilizer, and we shall assume that $U^c_h|_{G_n}=
U^{\tilde{c}}_{\tilde{h}}|_{G_n}$ (so in particular we assume that
the two representations of $\diff$ are given on the same Hilbert
space). To simplify notations, for the rest of the proof we shall
further set $U\equiv U^c_h$ and $\tilde{U}\equiv
U^{\tilde{c}}_{\tilde{h}}$.

Suppose $\xi \in \diff$ is such that it preserves all but one 
of our $n$ fixed points. Let this point be $p_j$. Set $q\equiv
\xi(p_j)$, and let us think of $G_{n+1}$ as the stabilizer of the
points $p_1,...,p_n$ and $q$; then $\xi^{-1} G_{n+1}\xi\subset
G_{n}$. Accordingly, we have that for all $\varphi\in G_{n+1}$
\begin{equation}
{\rm
Ad}\left(U(\xi)\tilde{U}(\xi^{-1})\right)(U(\varphi))=U(\varphi).
\end{equation}
However, we cannot immediately conclude that
$U(\xi)\tilde{U}(\xi^{-1})$ commutes with $U(\varphi)$, since the
above equation is meant in the sense of {\it projective}
unitaries. Nevertheless, it follows that there exists a complex
unit number $\lambda(\xi)$, such that in the sense of unitary
operators (i.e.\! not only in the projective sense) we have
\begin{equation}
U(\varphi)^* \,{\rm
Ad}\left(U(\xi)\tilde{U}(\xi^{-1})\right)(U(\varphi)) =
\lambda(\xi)\mathbbm 1.
\end{equation}
Clearly, the value of $\lambda(\xi)$ is independent of the chosen
phase of $U(\varphi)$, and moreover, it is largely independent
from the diffeomorphism $\xi$. Indeed, if $\xi'$ is another
diffeomorphism such that $\xi'(p_k)=p_k$ for $k\neq j$ and
$\xi'(p_j)=q$, then $\xi'=\xi\circ\beta$ where $\beta\equiv
\xi^{-1}\circ\xi' \in G_n$ and thus $U(\beta)=\tilde{U}(\beta)$
and so in the projective sense
\begin{equation}
U(\xi')\tilde{U}(\xi'^{-1}) =U(\xi)\tilde{U}(\beta)\,
\tilde{U}(\xi'^{-1}) = U(\xi)\tilde{U}(\xi^{-1})
\end{equation}
implying that $\lambda(\xi')=\lambda(\xi)$. However, the map $\xi
\mapsto \lambda(\xi)$ is clearly continuous, so the above argument
actually shows that $\lambda(\xi)=1$; i.e.\! that
$U(\xi)\tilde{U}(\xi^{-1})$ commutes with $U(\varphi)$. Hence we
have shown that $U(\xi)\tilde{U}(\xi^{-1})$ is in the commutant of
$U(G_n)$ and so --- by the condition of irreducibility --- it
follows that $U(\xi)=\tilde{U}(\xi)$. This concludes our proof,
since $\diff$ is evidently generated by the diffeomorphisms that
preserve all but one of the points $p_1,...,p_n$.
\end{proof}
At this point it is natural to ask: what are the admissible pairs
$(c,h)$, for which we can prove the irreducibility of
$U^c_h|_{G_n}$ for some $n>1$? (Recall that for $n=1$ we have
already obtained irreducibility, but in order to use the above
proposition, we need $n>1$.)

Here we shall prove irreducibility for two (overlapping) regions:
for $c\leq 1$, and for $h=0$ (the latter only for $n\leq 3$). The
irreducibility in the first of them is an evident consequence of
the known properties of the Virasoro nets. Nevertheless, it is
worth to state it.
\begin{proposition}
\label{prop:c<1} Let $(c,h)$ be an admissible pair with $c\leq 1$.
Then $U^c_h|_{G_n}$ is irreducible for every positive integer $n$.
Moreover, $U^c_h|_{G_n}\simeq U^c_{\tilde{h}}|_{G_n}$ if and only
if $h=\tilde{h}$.
\end{proposition}
\begin{proof}
As it is known, \cite{kl,xu03}, the Virasoro net with $c\leq 1$ is
strongly additive. Moreover, it is also known, if $(c,h)$ is an
admissible pair with $c\leq 1$, then the representation $U^c_h$
gives rise to a {\it locally normal} representation of the
conformal net $\A_{{\rm Vir}_c}$; see the discussion before
\cite[Prop.\! 2.1]{Car} explaining for which values of $c$ and $h$
it is known (and from where) that $U^c_h$ gives a locally normal
representation of $\A_{{\rm Vir}_c}$. This clearly shows that
$U^c_h|_{G_n}$ is irreducible. The rest of the proposition follows
from irreducibility and the previous proposition.
\end{proof}
\begin{proposition}\label{prop:h=0}
Let $(c,h=0)$ be an admissible pair and $n\leq 3$. Then
$U^c_0|_{G_n}$ is irreducible.
\end{proposition}
\begin{proof}
It is enough to show the statement for $n=3$. As usual in case of
$h=0$, we shall omit the index of the lowest energy, and we shall
denote the lowest energy vector by $\Omega$ (the ``vacuum
vector'') instead of $\Psi_0$. Moreover, we shall set $L_k\equiv
\pi^c_0(l_k)\; (k\in \ZZ)$.

The proof relies on the simple fact that in case of $h=0$, the
equality $L_k\Omega= L_k^+\Omega = 0$ is satisfied for $3$
different values of $k$; namely for $k=0,\pm 1$. Using this, we
shall show that each energy level of $V^c_0$ is in
\begin{equation}
S\equiv {\rm Span}\{A_1^+...A_j^+\Omega|j\in\NN, A_1,...A_j\in
\pi^c_0({\mathfrak K}_3),A_1\Omega=...=A_j\Omega=0\}
\end{equation}
(where $j=0$ means the vector $\Omega$ itself). This is enough;
then the statement follows by Prop.\! \ref{prop:irred}.

We shall argue by induction on the energy level. The zero energy
level $(V^c_0)_{(0)}$ is in $S$, since $(V^c_0)_{(0)}=\CC\Omega$.
So suppose that $(V^c_0)_{(k)}\subset S$ for all $k\leq m$, and
consider the case $k=m+1$.

Of course, the energy level $(V^c_h)_{(m+1)}$ is spanned by
vectors of the form $L_{-n_1}...L_{-n_j}\Omega$ where $j$ and
$n_1\leq  n_2...\leq n_j$ are positive integers such that
$n_1+...+n_j=m+1$. However, as it is well known, for $h=0$, these
vectors are not independent, and $(V^c_0)_{(m+1)}$ is already
spanned by the vectors of the above form with the further
condition that $2\leq n_1 \leq n_2...\leq n_j$.

So consider one of these vectors, and let $r$ be the number in
$\{0,\pm 1\}$ such that $r\equiv n_1$ modulo $3$. Then setting
$A\equiv L_{n_1}-L_r$ we have that $A\in \pi^c_0({\mathfrak K}_3)$
and $A\Omega=0$. It follows that $A^+S\subset S$. Moreover,
\begin{equation}
L_{-n_1}...L_{-n_j}\Omega = A^+(L_{-n_2}...L_{-n_j}\Omega) +
L_{-r}(L_{-n_2}...L_{-n_j}\Omega)
\end{equation}
and of course by the inductive condition both the vector
$L_{-n_2}...L_{-n_j}\Omega$ and the vector
$L_{-r}L_{-n_2}...L_{-n_j}\Omega$ is in $S$ (as $r<2\leq n_1$, the
energies of both vectors are smaller than $m+1$). Thus by the
above equation $L_{-n_1}...L_{-n_j}\Omega\in S$ and so
$(V^c_h)_{(m+1)}\subset S$, which concludes the inductive argument
and our proof.
\end{proof}

\section{Constructing representations of ${\mathfrak K}$}

Recall that $\theta(k_n)=k_{-n}$ and $\theta(C)=C$ and hence $\theta({\mathfrak K})=
{\mathfrak K}$. A representation $\eta$ of ${\mathfrak K}$ on complex
scalar product space $V$ satisfying $\eta(\theta(x))=\eta(x)^+$ for every $x\in
{\mathfrak K}$ will be said to be {\bf unitary}. So far, as a concrete example for
such representation, we only had the representations $\pi^c_h|_{\mathfrak K}$ obtained
by restriction. We shall now exhibit more examples.

Let us begin now our list of constructions with an abstract one.
Suppose that we have a $\gamma :{\mathfrak  K} \rightarrow
{\mathfrak K}$ endomorphism that commutes with the antilinear
involution $\theta$. Then it is clear, that for any unitary
representation $\eta$, the composition $\eta\circ\gamma$ is still
a unitary representation.

For example, following the similar constructions for the Virasoro
algebra, cf.\! \cite{LX}, for any $r\in \NN^+$ consider the linear
map $\gamma_r$ given by
\begin{eqnarray}
\label{endomorphism}
\nonumber
k_n &\mapsto& \frac{1}{r} k_{rn} + \frac{C}{24}(r-\frac{1}{r}) \;\;\;
(n\in\ZZ\setminus\{0\}) \\
C &\mapsto& r C.
\end{eqnarray}
By a straightforward calculation using the commutation relations
of the algebra $\mathfrak K$, we have that $\gamma_r$ is an
endomorphism and it is clear that it commutes with $\theta$. Thus
for any unitary representation we can construct a family of new
unitary representations by taking compositions with $\gamma_r$.

Just as in the case of the Virasoro algebra, we can also get some interesting
constructions considering the $U(1)$ current algebra. As it is well-known,
for every $q\in \RR$ there exists a linear space $V_q$ with positive scalar
product, a unit vector $\Phi_q\in V_q$ and a set of operators
$\{J_n \in {\rm End}(V_q)| n\in \ZZ\}$ satisfying the following
properties.
\begin{itemize}
\item
$[J_n,J_m] = n \,\delta\!_{-n,m} \,\mathbbm 1$ and
$J_{-n}=J_n^+$.
\item
$J_n \Phi_q = 0$ for all $n>0$.
\item
$J_0 = q \mathbbm 1$.
\item
$V_q$ is the smallest invariant subspace for $\{J_n|n\in\ZZ\}$,
containing $\Phi_q$.
\end{itemize}
We shall call this representation of the $U(1)$ current algebra
the representation with charge $q\in\RR$. The formally infinite
sum of the {\it normal product} of the current with itself
\begin{equation}
:\!J^2\!:_n \equiv \sum_{k>n-k} J_{n-k}J_k + \sum_{k\leq n-k}J_k J_{n-k}
\end{equation}
becomes finite on each vector of $V_q$, thus giving a well-defined
linear operator. Setting $L_n \equiv \frac{1}{2}:\!J^2\!:_n$ one finds
that the map $l_n\mapsto L_n$ extends to a unitary representation
of the Virasoro algebra with the central charge represented by
$\mathbbm 1$. Moreover, one finds that
\begin{itemize}
\item
$[L_n,J_m] = -m \,J_{n+m}$,
\item
$L_n\Phi_q = 0$ for all $n>0$,
\item
$L_0\Phi_q = \frac{1}{2}\,q^2 \,\Phi_q$.
\end{itemize}
We shall now give a new construction for some unitary
representation of $\mathfrak K$. The next proposition --- although
it can be understood and justified even without knowing anything
more than what was so far listed --- needs some ``explanations''.
Without making explicit definitions and rigourous arguments, let
us mention the following. (In any case, the precise statement and
its proof will make no explicit use of this.)

The main idea of \cite{BS-M} is the fact, that --- using their
settings\footnote{After a Cayley transformation, the formula
on the real lines simplifies to ``stress-energy tensor + $\alpha$-times the
(real-line) derivative of the current''.} and notations, but changing
the singular point from $-1$ to $1$ --- on the punctured
plane $\CC\setminus\{1\}$,
\begin{equation}
\label{formulaBS-M}
T_\alpha(z) = T(z) + \alpha \left(J'(z) + i\frac{z+1}{z-1}J(z)\right)
\end{equation}
follows the commutation relations of a stress-energy tensor at central charge
$c=1+12\alpha^2$. However, the function $z\mapsto i\frac{z+1}{z-1}$ has a
singularity at point $z=1$, and the power series expansion of $T_\alpha$
will depend on the chosen region. As a consequence, the operators
appearing in the expansion will give rise to a representation of the
Virasoro algebra, which --- on the full representation space $V_0$ ---
will not satisfy the unitarity condition. On the other hand,
with $h(z) \equiv 1-z^n$, the product $h T_\alpha$ will have an
unambiguous expansion, as in fact
\begin{equation}
\frac{z+1}{z-1} h(z) =
- (z+1)(1+ z + z^2 \ldots z^{n-1}) = 1 + z^n - 2 \sum_{k=0}^{n} z^k
\end{equation}
is a polynomial. This suggests the following statement.
\begin{proposition}
\label{mainconstruction}
For any fixed $\alpha\in \RR$, setting
$$
K^{\alpha}_n \equiv (L_0-L_n) + i n \alpha \left(J_n
+ \frac{1}{|n|}(J_0 +J_n - 2\!\!\!\!\!\sum_{k={\rm min}(0,n)}^{{\rm
max}(0,n)}\!\!\!\!\! J_k\,)\right)
$$
the assignment $k_n\mapsto K^{\alpha}_n\;(n\in\ZZ\setminus \{0\}), \, C\mapsto
c(\alpha)=1 + 12 \alpha^2$ gives a unitary representation of $\mathfrak K$.
\end{proposition}
Unitarity is manifest, and the rest of the proposition may be justified
by a long, but straightforward calculation using what was previously
listed about the $U(1)$ current. Note that here the formula is given in a
compact way, which is a fine thing for a proposition, but not necessarily
the best for actual calculations.  (For example, there is a hidden sign
factor, appearing as $n/|n|$.)

The reader is encouraged to {\it really} check the commutation
relations. It might seem something tedious (and boring), but ---
according to the author's personal opinion --- in fact it is
interesting to observe how the apparent contradictions disappear
by some ``miraculous cancellations'' of the terms.

\section{Equivalence with $\mathbf{h_1\neq h_2}$}

Let $(c,h)$ be an admissible pair and consider the representation
$\pi^c_h$ on the representation space $V^c_h$. It is clear, that
for any $\Phi\in V^c_h$ we have that
$\pi^c_h(k_n)\Phi=(\pi^c_h(l_0)-\pi^c_h(l_n))\Phi$ which, for $n$
sufficiently large, is further equal to $\pi^c_h(l_0)\Phi$. This
shows two things: first, that up to phase $\Phi=\Psi^c_h$ is the
unique normalized vector in $V^c_h$ such that $\pi^c_h(k_n)\Phi =
h \Phi$ for all $n>0$; second, that by knowing the restriction
$\pi^c_h|_{\mathfrak K}$ one can ``recover'' the operator
$\pi^c_h(l_0)$ and hence the whole representation $\pi^c_h$. It
follows that $\pi^c_h|_{\mathfrak K}$ is equivalent to
$\pi^{\tilde{c}}_{\tilde{h}}|_{\mathfrak K}$ if and only if
$\pi^c_h$ is equivalent to $\pi^{\tilde{c}}_{\tilde{h}}$; i.e.\!
if and only if $(c,h)=(\tilde{c},\tilde{h})$. However, as it was
already explained, what we are interested in is not the
equivalence of $\pi^c_h|_{\mathfrak K}$ and
$\pi^{\tilde{c}}_{\tilde{h}}|_{\mathfrak K}$ but the equivalence
of $\overline{\pi^c_h}|_{\mathfrak K}$ and
$\overline{\pi^{\tilde{c}}_{\tilde{h}}}|_{\mathfrak K}$. In order
to investigate the second kind of equivalence, first we shall
consider different ways of exhibiting the representation
$\pi^c_h|_{\mathfrak K}$.

Let $\eta$ be a representation of $\mathfrak K$ on a complex scalar
product space $V$. Recall that $\theta(k_n)=k_{-n}, \theta(C)=C$ and so
$\theta({\mathfrak K})={\mathfrak K}$. Assume that $\eta$ satisfies the
following properties:

\begin{itemize}

\item[(A)] $\eta$ is unitary:
$\eta(\theta(x))=\eta(x)^+$ for all $x\in {\mathfrak K}$,

\item[(B)]
$\eta(C) = c \mathbbm 1$,

\item[(C)] up to phase there exists a unique
normalized vector $\Psi$ with the property $\eta(k_n)\Psi = h \Psi$ for
all $n>0$,

\item[(D)] $V$ is the smallest invariant space for $\eta$
containing $\Psi$.
\end{itemize}
Using the commutation relations and the listed properties it is an
exercise to show that

\begin{itemize}

\item The value of the scalar product
$$\langle \eta(k_{n_1})\eta(k_{n_2})\ldots \eta(k_{n_r}) \Psi,\,
\eta(k_{m_1})\eta(k_{m_2})\ldots \eta(k_{m_s}) \Psi\rangle$$ is
``universal'': it is completely determined by the values of $c,h$
and the integers $n_1,\ldots ,n_r$ and $m_1,\ldots ,m_s$. That is,
the scalar product can be calculated by knowing these values; even
without having the actual form of the representation $\eta$ or
knowing anything more (than just the required properties) about
it.

\item The representation space is
spanned by the vectors of the form appearing in Lemma
\ref{vermalike-module}.
\end{itemize}
These two consequences imply that the representation, up to equivalence,
is uniquely determined by the pair $(c,h)$. It is worth to state this in a
form of a statement.
\begin{corollary} \label{uniqlowesterep}
Let everything be as it was
explained. Then the map
$$\eta(k_{n_1})\eta(k_{n_2})\ldots\eta(k_{n_r})\Psi \,\mapsto \,
\pi^c_h(k_{n_1})\pi^c_h(k_{n_2})\ldots\pi^c_h(k_{n_r})\Psi^c_h $$
extends to a unique unitary operator which establishes an
isomorphism between $\eta$ and $\pi^c_h|_{\mathfrak K}$.
\end{corollary}
We shall now get to the ``main trick'' of this paper, which is a
combination of the two constructions discussed in the previous section. So
consider the unitary representation of $\mathfrak K$ given by Proposition
\ref{mainconstruction} for $q=0$, and compose it with the endomorphism
$\gamma_2$ given by (\ref{endomorphism}). We get that for every
$\alpha \in\RR$, the map $k_n\mapsto K^{(\alpha,2)}_n \,
(n\in(\ZZ\setminus \{0\}))$, where
\begin{eqnarray} \nonumber
K^{(\alpha,2)}_n \,\equiv & &\frac{1}{2}(L_0-L_{rn}) \\ \nonumber & +& i
n\alpha \left(J_{2n} + \frac{1}{2|n|}(J_{2n} - 2\!\!\!\!\!\!\sum_{k={\rm
min}(0,2n)}^{{\rm max}(0,2n)}\!\!\!\!\! J_k\,)\right) \\ &+&
\frac{1+12\alpha^2}{16}\mathbbm 1
\end{eqnarray}
extends to a unitary representation $\rho_{(\alpha,2)}$ of $\mathfrak K$
with central charge $c(\alpha,2) = 2 (1+12\alpha^2)$ (i.e.\, the element
$C$ is represented by $c(\alpha,2)\mathbbm 1 $). Note that in the above
formula we omitted $J_0$, since we are in the {\bf
vacuum representation} of the $U(1)$ current (i.e.\! $q=0$ and so $J_0=0$,
too). Moreover, as it is usual in the vacuum representation, we shall
denote denote the lowest energy vector --- corresponding to the ``true
conformal energy'' $L_0$ --- by $\Omega$, rather than by $\Phi_0$, and we
shall call it the {\bf vacuum vector}.

\begin{lemma}
Let $\Phi\equiv J_{-1}\Omega + i \frac{16\alpha}{1+12\alpha^2} \Omega$.
Then for every $n>0$ we have
$$K^{(\alpha,2)}_n \Omega \,= \,\frac{1+12\alpha^2}{16}\,\Omega,
\;\;\;\;
K^{(\alpha,2)}_n \Phi \,= \,\frac{9+12\alpha^2}{16}\,\Phi.$$
\end{lemma}
\begin{proof}
Anything which lowers the energy (in the sense of $L_0$) by more
than $1$, annihilates both the vector $\Omega$ and $\Phi$. Thus
$J_k\Phi = L_k\Phi=J_k\Omega=L_k\Omega = 0$ for every $k>1$ and
hence one finds that the operator $K^{(\alpha,2)}_n$, for every
$n>0$, acts exactly like the operator
\begin{equation}
\frac{1}{2}L_0 - i\alpha J_1 + \frac{1+12\alpha^2}{16}\mathbbm 1
\end{equation}
on the mentioned vectors. The rest is trivial calculation.
\end{proof}

By the previous lemma and by Corollary \ref{uniqlowesterep}, for
$c>2, h_1 = \frac{c}{32}$ and $h_2 = \frac{1}{2}+\frac{c}{32}$,
the representations $\pi^c_{h_1}|_{\mathfrak K}$ and
$\pi^c_{h_2}|_{\mathfrak K}$ appear as subrepresentations of a
common, non-irreducible representation, namely the representation
$\rho_{(\alpha,2)}$ with $\alpha = \sqrt{\frac{c/2-1}{12}}$.

Let $V_{h_1}$ be the minimal invariant subspace for $\rho_{(\alpha,2)}$
containing $\Omega$, and $V_{h_2}$ the one containing the previously given
vector $\Phi$. These are the subspaces on which the representation is
isomorphic to $\pi^c_{h_1}|_{\mathfrak K}$ and $\pi^c_{h_2}|_{\mathfrak
K}$, respectively, since --- as we have seen --- the vectors $\Omega$
and $\Phi$ behave like ``lowest energy vectors'' for the representation
$\rho_{(\alpha,2)}$ of $\mathfrak K$.

The important observation is that these two vectors, since $\alpha\neq 0$,
are not orthogonal. Thus, neither the two subspaces $V_{h_1}$ and $V_{h_2}$
can be so. It {\it should} follow therefore, that the corresponding
irreducible representations (consisting of closed operators) cannot be
inequivalent.

This argument however, is not completely rigourous as we deal with
(unbounded) operators rather than unitary representations of groups. So in
what follows, we shall find a way to deal with the technical difficulties.

\begin{lemma}
Let $Q_1$ and $Q_2$ be the orthogonal projections onto
$\overline{V}_{h_1}$ and $\overline{V}_{h_2}$. Then for every
$x\in {\mathfrak K}$ such that $\theta(x)=x$, we have that
$\overline{\rho_{(\alpha,2)}(x)}$ is self-adjoint and
$$
Q_j e^{ i t \,\overline{\rho_{(\alpha,2)}(x)}  }=
e^{it\,\overline{\rho_{(\alpha,2)}(x) } } Q_j \;\;\;(t\in \RR, j=1,2).
$$
\end{lemma}
\begin{proof}
The finite energy vectors of the $U(1)$ current are analytic for every
operator which is a finite sum of the ``$L$'' operators; in particular,
for $\rho_{(\alpha,2)}(x)$, too. This shows, that if $\theta(x)=x$ then
$\rho_{(\alpha,2)}(x)$ is essentially self-adjoint. Moreover, as the
the subspaces $V_{h_j}\;(j=1,2)$ are invariant for
$\rho_{(\alpha,2)}(x)$, the analyticity property also shows, that
the subspaces $\overline{V}_{h_j}\;(j=1,2)$ are invariant for
$e^{it\,\overline{\rho_{(\alpha,2)}(x)}}\;(t\in \RR)$.
\end{proof}

\begin{corollary}
Let $c>2, h_1 = \frac{c}{32}$ and $h_2 = \frac{1}{2}+\frac{c}{32}$. Then
the representations $U^c_{h_1}|_{G_1}$ and  $U^c_{h_2}|_{G_1}$
are unitary equivalent.
\end{corollary}

\begin{proof}
By Corollary \ref{algebra->group}, the restrictions of
$\rho_{(\alpha,2)}$ onto $V_{h_1}$ and $V_{h_2}$ give rise to two
unitary representations (in that corollary, all elements of
$\mathfrak K$ appear in the condition, but it is clear, that the
hermitian ones, i.e.\! the elements invariant under $\theta$, are
sufficient for us) of $G_1$ on $\overline{V}_{h_1}$ and
$\overline{V}_{h_1}$, respectively, with the first one being
unitary equivalent to $U^c_{h_1}|_{G_1}$, while the second one to
$U^c_{h_2}|_{G_1}$. By Corollary \ref{cor:irred}, these
representations are irreducible, and by the previous lemma and
Prop.\! \ref{intertwiners}, the restriction of $Q_1$ is an
intertwiner between them. Hence if $Q_1$, as a map from
$\overline{V}_{h_2}$ to $\overline{V}_{h_1}$, is not zero, then
the two representations are equivalent. This is indeed the case,
as $\langle \Omega,\Phi\rangle = \frac{16i\alpha}{1+12\alpha^2}
\neq 0$, and $\Omega\in V_{h_1}$ whereas $\Phi\in V_{h_2}$.
\end{proof}

Thus we have managed to give examples for values $h\neq\tilde{h}$
such that $U^c_h|_{G_1}\simeq U^{c}_{\tilde{h}}|_{G_1}$. More
examples could be generated by $i$) taking tensor products (and
then restrictions), $ii$) using the endomorphism $\gamma_r$ with
$r$ different from $2$ (which we have used so far). However, at
the moment our aim was just to find {\it some} examples.
\begin{corollary}\label{cor:localglobal}
In the Virasoro net with $c>2$, local and global intertwiners are
not equivalent.
\end{corollary}
\begin{proof}
The representations $U^c_{h}$ and $U^c_{\tilde{h}}$, where $h =
\frac{c}{32}$ and $\tilde{h} = \frac{1}{2}+\frac{c}{32}$ give rise
to two locally normal, irreducible representations of the
conformal net $\A_{{\rm Vir}_c}$; see the discussion before
\cite[Prop.\! 2.1]{Car} explaining for which values of $c$ and $h$
it is known (and from where) that $U^c_h$ gives a locally normal
representation of $\A_{{\rm Vir}_c}$.

In any locally normal irreducible representation of $\A_{{\rm
Vir}_c}$, there is a unique strongly continuous unitary
representation of the universal covering group of the M\"obius
group, which implements the M\"obius symmetry in the given locally
normal irreducible representation of the net; see
\cite{D'AnFreKos} for the details. This shows that the value of
the lowest energy in any locally normal irreducible representation
of $\A_{{\rm Vir}_c}$ is well-determined and hence, {\it
globally}, the locally normal irreducible representations given by
$U^c_{h}$ and $U^c_{\tilde{h}}$ are inequivalent. However, by the
previous result, {\it locally} they are equivalent.
\end{proof}

\noindent {\bf Acknowledgements.} The author would like to thank
Roberto Longo for suggesting the problem. The main construction of
this paper, providing a counter-example for the equivalence of
local and global intertwiners, was found in March 2006, while the
author stayed at the ``Institute f\"ur Theoretische Physik'' in
G\"ottingen (Germany) and was supported in part by the EU network
program ``Quantum Spaces -- Noncommutative Geometry'' and by the
``Deutsche Forschungsgemeinschaft''. The author would like to
thank the hospitality of the Institute and in particular the
hospitality of Karl-Henning Rehren and moreover, the useful
discussions with him that contributed in a significant way to
finding the main construction. The author would also like to thank
Sebastiano Carpi for further useful discussions. Finally, the
author would like to thank the organizers of the conference
``Recent Advances in Operator Algebras'' held in Rome, November
8-11, 2006 (on the occasion of the 60th birthday of Laszl\'o
Zsid\'o), where the result presented here was first announced.

\end{document}